\documentclass{amsart}
\usepackage{amssymb}
\newcommand{\Miller}{{\mathbb M}}
\newcommand{\Cohen}{{\mathbb C}}

\newcommand{\stem}{{\operatorname{\mathsf   {stem}}}}
\newcommand{\spli}{{\operatorname{\mathsf   {split}}}}
\newcommand{\forces}{\Vdash}
\newcommand{\suc}{{\operatorname{\mathsf    {succ}}}}
\newcommand{\ZFCa}{{\operatorname{\mathsf {ZFC}}}}
\newcommand{\CH}{\operatorname{\mathsf {CH}}}
\newcommand{\reals}{{\mathbb R}}

\newcommand{\rest}{{\mathord{\restriction}}}
\newcommand{\co}{{\mathbf c}}

\newcommand{\ran}{{\operatorname{\mathsf    {range}}}}

\newcommand{\V}{{\mathbf V}}
\newcommand{\<}{\langle}
\renewcommand{\>}{\rangle}
\newcommand{\thinks}{\models}
\newcommand{\xp}{{\mathbf X}}

\newcommand{\lft}[2]{\mathopen\ifcase#1{}\oo\or
                        \big#2\or\Big#2\else\oo\fi} 
\newcommand{\rgt}[2]{\mathclose\ifcase#1{}\oo\or
                        \big#2\or\Big#2\else\oo\fi} 
\newcommand{\PM}{{\mathbf {PM}}}
\newcommand{\UM}{{\mathbf  {UM}}}
\newcommand{\bor}{{\mathsf   {Borel}}}
\newcommand{\Pp}{{\mathsf   {P}}}

\theoremstyle{plain}
\newtheorem{theorem}{Theorem}
\theoremstyle{plain}
\newtheorem{lemma}[theorem]{Lemma}

\newtheorem{definition}[theorem]{Definition}

\begin{document}
\title{On perfectly meager sets}
\author{Tomek Bartoszynski}
\address{Department of Mathematics and Computer Science\\
Boise State University\\
Boise, Idaho 83725 U.S.A.}
\thanks{Author  was partially supported by 
NSF grant DMS 99271}
\email{tomek@math.idbsu.edu, http://math.idbsu.edu/\char 126 tomek}

\begin{abstract}
We show that it is consistent that 
the product of perfectly meager sets is perfectly meager.
\end{abstract}
\maketitle

\section{Introduction}

Suppose that $X$ is a subset of a Polish space $\xp$. We say that   
$X$ is {\em perfectly meager} if for every perfect set $P
  \subseteq \xp $, $X \cap P$ is meager in the relative topology of
  $P$. Let $\PM$ denote the collection of perfectly meager sets.

Clearly all countable sets are perfectly  meager but there are various
examples of uncountable perfectly meager sets that can be constructed
in $\ZFCa$ (see \cite{Mil84Spe}).

In \cite{Mar68}, Marczewski asked whether the product of perfectly meager
sets is perfectly meager.
This question was partially answered by Rec{\l}aw who showed that:
\begin{theorem}[Rec{\l}aw \cite{Rec91Pro}]
  Assume $\CH$. Then there are two perfectly meager sets whose product
  is not perfectly meager.
\end{theorem}

The proof relies on the existence of a Borel set having certain
properties and 
the existence of a Luzin
set (i.e. an uncountable set whose intersection with every meager set
is countable).

The purpose of this note is to show that it is also consistent with
$\ZFCa$ that the product of any two perfectly meager sets is perfectly
meager.
Thus, Marczewski's question is undecidable in $\ZFCa$.

\begin{definition}
  Suppose that $X$ is a subset of a Polish space $\xp$. We say that
$X$ is {\em universally meager} if every Borel isomorphic image
  of $X$ in $\xp$ is meager.
Let $\UM$ denote the
collection of universally meager sets.
\end{definition}

It is clear that $\UM \subseteq \PM$, but the other inclusion may fail.
\begin{theorem}[Sierpi\'nski \cite{Sie}]
  Assume $\CH$. Then $\PM \neq \UM$.
\end{theorem}

Unlike $\PM$, the class of universally meager sets is closed under
products.

\begin{theorem}[Zakrzewski \cite{Zakr}]
  The product of universally meager sets is universally meager.
\end{theorem}

This is a consequence of the following characterization of the class
$\UM$.
Let $\Cohen$ denote the Cohen algebra.
\begin{theorem}[Zakrzewski \cite{Zakr}] 
For a subset $X$ of a perfect Polish space $\xp$, the
following are equivalent:
\begin{enumerate}
\item $X \in \UM$.
\item $X$ does not contain a Borel one-to-one image of a non-meager set.
\item For every $\sigma$-ideal $ {\mathcal J}$ in $\bor(\xp)$ such that
  $\bor(\xp)/ {\mathcal J} 
\cong \Cohen$ there is a Borel set $B\in {\mathcal J} $ such that $X\subseteq
B$.
\item $X$ is meager in every  Polish topology $\tau$
on $\xp$ such that $\xp$ has no isolated points and $\bor(\xp,\tau)=\bor(\xp)$.
\item $X$ is meager in every  second countable Hausdorff topology
$\tau$ on $\xp$ such that $\xp$ has no isolated points and all Borel sets
(in the original Polish
topology) have Baire Property in the topology~$\tau$.
\item There is no $\sigma $-ideal ${\mathcal J}$ in $\bor(X)$ such that
$\bor(X)/{\mathcal J} \cong \Cohen$.
\item $X$ is meager in every second countable Hausdorff topology
$\tau$ on $X$ such that $X$ has no isolated points and all Borel
subsets of $X$ (in the topology
inherited from the original Polish
topology on $\xp$) have Baire Property in the topology~$\tau$.
\item $X$ is meager in every separable metrizable topology
$\tau$ on $X$ such that $X$ has no isolated points and
$\bor(X)=\bor(X,\tau)$.
\end{enumerate}
 \end{theorem}

Now we can formulate our main result:
\begin{theorem}
  It is consistent that $\PM=\UM$. 
In particular, it is consistent that $\PM$ is closed under products.
\end{theorem}

\section{In $\ZFCa$}
In this section we will identify a more general property that implies
that $\PM=\UM$.

For a function $f: X \longrightarrow Y$ let $f[X]=\{f(x):x\in X\}$
denote the image of $X$.

We need the following observation:
\begin{lemma}\label{one}
  Suppose that $X \not \in \UM$. Then there exists  $X' \subseteq X$,
  a Borel set $A \supseteq X'$ and Borel isomorphism $f:A
  \longrightarrow \omega^\omega $ such that 
  \begin{enumerate}
  \item $f^{-1}$ is continuous,
  \item $f[X']$ is not meager in $\omega^\omega $.
  \end{enumerate}
\end{lemma}
\begin{proof}
Without loss of generality we can assume that $X \subseteq \reals=\xp$.
Let $h:X \longrightarrow h[X]$ be the Borel isomorphism witnessing
that $X \not \in \UM$.
Find  Borel sets $A \supseteq X$ and $B \supseteq h[X]$ and the Borel
isomorphism $\bar{h}:A \longrightarrow B$ extending $h$ (exercise in
\cite{Srivastava}). 
By removing a countable set from $A$ we can arrange that $B$ is $0$-dimensional. 
Let $F \subseteq B$ be a meager set such that 
\begin{enumerate}
\item $\bar{h}[X] \cap (B \setminus F) $ is not meager in $B \setminus F$,
\item $\bar{h}^{-1} \rest (B \setminus F) $ is continuous on $B \setminus F$,
\item $B \setminus F$ is homeomorphic to $\omega^\omega $.
\end{enumerate}

Apply Kuratowski's theorem to get (2). Without loss of generality we
can assume that $F$ is an $F_\sigma $ set which is dense.
By the theorem of Mazurkiewicz \cite{Kur66Top}, a $G_\delta $ subset
of a $0$-dimensional Polish space which has  empty interior is
homeomorphic to $\omega^\omega $. Call this homeomorphism $g$. Now the
required mapping is $g \circ \bar{h}$.
\end{proof}

Consider the following principle:

\begin{definition}
Axiom $\Pp$: For every nonmeager subset  $X \subseteq \omega^\omega $ there exists a
compact subset $P \subseteq \omega^\omega $ such that $P \cap X$ is
nonmeager in $P$.
\end{definition}
\begin{theorem}
  Assume $\Pp$. Then $\UM=\PM$.
\end{theorem}
\begin{proof}
  Suppose that $X \not\in \UM$.
By lemma \ref{one} there exists $X' \subseteq X$, 
a Borel set $A \supseteq X'$ and Borel isomorphism $f:A
  \longrightarrow \omega^\omega $ such that 
  \begin{enumerate}
  \item $f^{-1}$ is continuous,
  \item $f[X']$ is not meager in $\omega^\omega $.
\end{enumerate}
Apply Axiom $\Pp$ to to find a compact set $P \subseteq \omega^\omega $
such that $P \cap f[X']$ is not meager in $P$.
Set $Q=f^{-1}(P)$ and note that $f \rest Q$ is a homeomorphism between
$P$ and $Q$. Under this homeomorphism 
$Q \cap X'$ is the image of $f[X'] \cap P$ thus it is
not meager in $Q$. Since $Q \cap X \supseteq Q \cap X'$ it follows
that $X \not\in \PM$. As the other inclusion is obvious, the theorem
follows.
\end{proof}

\section{Forcing}
In this section we will show that Axiom $\Pp$ is consistent with $\ZFCa$,
which will finish the proof.

For a tree $p$ and $t \in p$,
let $\suc_p(t)$ be the set of
all immediate successors of $t$ in $p$, $p_t=\{v \in p: t
\subseteq v \text{ or } v \subseteq t\}$ the subtree of $p$ determined
by $t$, and let $[p]$ be the set of
branches of $p$. By identifying $s \in \omega^{<\omega}$ with the
full-branching tree having root $s$, we can also denote
$[s]=\{f \in \omega^\omega:s \subseteq f\}$.
Let $\omega^{\uparrow \omega } = \{s \in \omega^{<\omega}: \text{$s$
  is strictly increasing}\}$.

The rational perfect forcing $\Miller$ is the following forcing notion:
\begin{multline*}
p \in \Miller \iff p \subseteq \omega^{\uparrow\omega} \text{ is a
  perfect tree } \& \\ 
\forall s \in p \ \exists t \in p \
\lft1(s \subseteq t \ \& \ |\suc_p(t)|=\boldsymbol\aleph_{0}\rgt1).
\end{multline*} 
For $p, q \in \Miller$, $p \geq q$ if $p \subseteq q$.
Without loss of generality we can assume that  $|\suc_p(s)|=1$ or
$|\suc_p(s)|=\boldsymbol\aleph_0$ for all $p\in \Miller$ and $s \in
p$.
Conditions of this type form a dense subset of $\Miller$.

Let      
$$ \spli(p) = \{s\in p : |\suc_p(s)|>1\} = \bigcup_{n \in \omega} \spli_n(p),$$
where
$\spli_n(p)=\left\{s \in \spli(p): \lft2|\lft1\{t \subsetneq s: t \in
\spli(p)\rgt1\}\rgt2| =n\right\}$.

For $p,q\in \Miller$, $n \in \omega$,   we let 
$$p \geq_n q \iff p \geq q \ \& \ \spli_n(q)=\spli_n(p).$$

If $v \in \spli(p)$ let $U^p_v=\{n \in \omega : v^\frown n \in p\}$
and for $n \in U^p_v$ let $v^n$ be the first splitting node below
$v^\frown n$.

If $G \subseteq \Miller$ is a generic filter over $\V$ let ${\mathbf
  m}=\bigcap_{p \in G} [p]$ be the generic real.

Let $\Miller_{\omega_2}$ be the countable support iteration of
$\Miller$ of length ${\boldsymbol\aleph}_2 $.

The following facts about $\Miller$ are well-known.
\begin{theorem}\label{wn}
  \begin{enumerate}
  \item The sequence $\< \leq_n: n \in \omega\>$ witnesses that
$\Miller$ satisfies axiom A. In particular, $\Miller$ is proper.
\item $\Miller$ preserves nonmeager sets, i.e. if $A \subseteq
  2^\omega $, $A \in \V$ is not meager then $\V^{\Miller} \thinks A
  \text{ is not meager}$. {\em (}\cite{BJbook}, theorem 7.3.46{\em )}
\item $\Miller$ also satisfies iterable condition for preserving
  non-meager sets. In particular, 
countable support iteration of Miller forcing preserves
  nonmeager sets {\em (}\cite{BJbook}, theorems 6.3.19 and 6.3.20{\em )}. \qed
  \end{enumerate}
\end{theorem}

\begin{theorem}\label{main}
  $\V^{\Miller_{\omega_2}} \thinks \text{Axiom $\Pp$}$.
\end{theorem}

The idea of the proof is as follows. Suppose that $X \in
\V^{\Miller_{\omega_2}}$  and $X$ is not meager in $\omega^\omega $. 
First we find $\alpha < \omega_2$ such that $\V^{\Miller_\alpha }
\thinks X \cap \V^{\Miller_\alpha } 
  \text{ is not meager}$ (lemma \ref{non}). Next we will find a compact set $P \subseteq
  \omega^\omega $ belonging to $\V^{\Miller_{\alpha+1} }$ such that 
$$\V^{\Miller_{\alpha+1} } \thinks P \cap X \cap \V^{\Miller_\alpha } 
  \text{ is not meager in  $P$ (theorem \ref{crucial})}.$$
Finally, by \ref{wn}(3),  $\Miller_{\omega_2}$ preserves
non-meager sets, thus
 $$\V^{\Miller_{\omega_2}} \thinks  X \cap \V^{\Miller_\alpha }
  \text{ is not meager in }P,$$
which implies that $\V^{\Miller_{\omega_2}} \thinks  X \text{ is not
  meager in }P.$

\begin{lemma}\label{non}
  Suppose that $X \in \V^{\Miller_{\omega_2}}, \ X \subseteq
  \omega^\omega $ and $X$ is a not meager in
  $\V^{\Miller_{\omega_2}}$. Then there exists $\alpha < \omega_2$
  such that 
$\V^{\Miller_\alpha } \thinks X \cap \V^{\Miller_\alpha }
  \text{ is not meager}.$
\end{lemma}
\begin{proof}
 Let $\<\alpha_\xi: \xi<\omega_1\> $ be a continuous increasing
  sequence such that 
$$\V^{\Miller_{\omega_2}} \thinks X \cap \V^{\Miller_{\alpha_{\xi+1}}}
  \text{ is not covered by any meager set from
    }\V^{\Miller_{\alpha_{\xi}}}.$$
By properness, $ \alpha =\sup_{\xi<\omega_1} \alpha_\xi$, has the
required property.
\end{proof}

\begin{theorem}\label{crucial}
  Suppose that $X \in \V$, $X \subseteq \omega^\omega $ is a non-meager
  set. There is a compact set $P \subseteq \omega^\omega $, $P \in
  \V^{\Miller}$ such that 
$\V^{\Miller} \thinks X \text{ is not
  meager in }P.$
\end{theorem}
\begin{proof}
For the sake of clarity we will break the proof into three lemmas.
The main idea of the proof is already present in 
\cite{JuSh478}.

Let $\dot{\mathbf m}$ be the canonical name for an $\Miller$-generic
real and 
let $\Cohen$ be the Cohen forcing represented as $\omega^{<\omega}$ with
$\dot{\mathbf c}$ being the canonical name for the Cohen real.

  For $x \in \omega^\omega $ let $P_x=\{z \in \omega^\omega : \forall
  n \ z(n) \leq x(n)\}$. Note that $P_x$ is a compact set in
  $\omega^\omega $.

For two sequences $s \in \omega^{<\omega},\ t \in \omega^{\uparrow \omega}$ we say that
$(s,t)$ is {\em good} if $|s| = |t|$ and $s(i) \leq t(i)$ for $i < |t|$.

\begin{lemma}\label{two}
  Suppose that $p \in \Miller$, $s \in \Cohen$,
  and $\lft1(s\rest |\stem(p)|,\stem(p)\rgt1)$ is good.
Then there is a $\Cohen$-name $\dot{q}$ for an element of $\Miller$
such that 
\begin{enumerate}
\item $s \forces_{\Cohen} \dot{q} \geq_0 p$,
\item $(s, \dot{q}) \forces_{\Cohen \star \dot{\Miller}} \dot{\mathbf c} \in
  P_{\dot{\mathbf m}}.$
\end{enumerate}
\end{lemma}
\begin{proof}
Let ${\mathbf c} \supseteq s$ be a Cohen real over $\V$.
Working in $\V[{\mathbf c}]$ define
$$q=\{v \in p: (\co \rest |v|,v) \text{ is good}\}.$$ 
It is enough to check that $q \in \Miller$. In fact, we will show that
if $v \in \spli(p)$ and $v \in q$ then $v \in \spli(q)$.
  
Suppose that $v \in \spli(p)$ and $v \in q$. In particular,
$(\co \rest |v|, v)$ is good.
Let $s'=\co \rest |v|$. 
For $k \in \omega $ let
$$D_k=\{t \in \Cohen: s' \subseteq t \ \&\ \exists n>k \ (t\rest |v^n|,v^n)
\text{ is good}\}.$$
We show that $D_k$ is dense in $\Cohen$ below $s'$.
Take any $s'' \geq s' $ and let $n \in U^p_v
\setminus \max\lft1(\ran(s''),k\rgt1)$. If $|s''|>|v^n|$ then put
$t=s''$. Otherwise,
let $t \geq s''$ be such that
$|t|=|v^n|$ and
$$t(i)=\left\{
  \begin{array}{ll}
s''(i)& \text{if } i < |s''|\\
0 & \text{if } |s''| \leq i < |v^n|
  \end{array}\right. .$$
It is clear that $t \in D_k$.
By genericity, we conclude that the set
$$\{n: (\co \rest |v^n|, v^n) \text{ is good}\}$$ is infinite in
  $\V[\co]$. In particular, $v \in \spli(q)$.
Since $\co$ was arbitrary, it finishes the proof.
\end{proof}

\begin{lemma}
   Suppose that $p \in \Miller$, $s \in \Cohen$,
  and $\lft1(s\rest |\stem(p)|,\stem(p)\rgt1)$ is good.
  Let $\dot{F}$ be an
  $\Miller$-name for a closed nowhere dense subset of $P_{\dot{\mathbf
      m}}$. 
There exists a $\Cohen$-name $\dot{q}$ for an element of $\Miller$
such that 
\begin{enumerate}
\item $s \forces_{\Cohen} \dot{q} \geq_0 p$,
\item $(s, \dot{q}) \forces_{\Cohen \star \dot{\Miller}} \dot{\mathbf c} \in
  P_{\dot{\mathbf m}},$
\item $(s, \dot{q}) \forces_{\Cohen \star \dot{\Miller}} \dot{\mathbf c}
  \not\in \dot{F}$.
\end{enumerate}
\end{lemma}
\begin{proof}
  Let ${\mathbf m} \in [p]$ be an $\Miller$-generic real over $\V$, and
  let $F$ be the interpretation of $\dot{F}$ using ${\mathbf m}$.
In $\V[{\mathbf m}]$ define sequences $\<s_n: n \in \omega\>\in
\lft1(\omega^{<\omega}\rgt1)^\omega$ such
that 
\begin{enumerate}
\item $\forall v \in \prod_{j<n} {\mathbf m}(j) \ [v^\frown s_n] \cap F = \emptyset$.
\item $ \forall j < |s_n| \ {\mathbf m}(n+j) \geq s_n(j)$,
\end{enumerate}
Since $F$ is nowhere dense 
this definition is correct.
Going back to $\V$ we conclude that there is an
$\Miller$-name  $\<\dot{s}_n: n \in \omega\>$ such that 
\begin{enumerate}
\item $p \forces_{\Miller} \forall v \in \prod_{j<n} \dot{\mathbf
    m}(j) \ [v^\frown \dot{s}_n] \cap \dot{F}=\emptyset$,
\item $p \forces_{\Miller} \forall n \ (\dot{s}_n, \dot{\mathbf m}
  \rest [n, n+|\dot{s}_n|) \text{ is good}$.
\end{enumerate}

For each $n \in U_{\stem(p)}^p$ find a condition $p_n \geq p$ such
that 
\begin{enumerate}
\item there is a sequence $s_n$ such that 
$p_n \forces_{\Miller} \dot{s}_n =s_n$.
\item $\stem(p_n)\supseteq \stem(p)^n$,
\item $|\stem(p_n)| \geq n+|s_n|$.
\end{enumerate}
Observe that by the choice of $\<\dot{s}_n: n \in \omega \>$ it
follows that
$$(s_n, \stem(p_n) \rest [n,n+|s_n|) \text{ is good}.$$
Let $\co \supseteq s$ be a Cohen  real over $\V$.
Working in $\V[\co]$ define
$$A=\left\{n \in \omega: \lft1(\co \rest |\stem(p_n)|, \stem(p_n)\rgt1) \text{ is
  good and } \co \rest \lft1[n, n+|s_n|\rgt1)=s_n\right\}.$$
We will show that $A$ is infinite in $\V[\co]$.

For $k \in \omega $ let 
\begin{multline*}D_k=\left\{t \in \Cohen : s \subseteq t \ \&\ \exists n>k \ \lft2(
  \lft1(t\rest |\stem(p_n)|,
    \stem(p_n)\rgt1) \text{ is good}\right. \ \& \\
\left.  t \rest \lft1[n, n+|s_n|\rgt1)=s_n\rgt2)\right\}.
\end{multline*}
We show that $D_k$ is dense in $\Cohen$ below $s$.
Suppose that $s' \geq s$ and let $\ell=\max\lft1(\ran(s'), |s'|,k\rgt1)$.
Pick
$n \in U^p_{\stem(p)}
\setminus \ell$ and define 
$t \geq s'$  such that
$|t|=|\stem(p_n)|$ and
$$t(i)=\left\{
  \begin{array}{ll}
s'(i)& \text{if } i < |s'|\\
0 & \text{if } |s'|\leq i<n\\
s_n & \text{if } n \leq i < n+|s_n|\\
0 & \text{if } n+|s_n| \leq i < |\stem(p_n)|
  \end{array}\right. .$$ 
 Note that 
by the properties of $\<\dot{s}_n: n \in \omega \>$ it follows that
$t \in D_k$.
By genericity, for every $k\in \omega $ there is $n \in \omega$ such that $\co \rest n \in
D_k$, which implies that $A$ is infinite.

Let $p^\star = \bigcup_{n \in U^p_{\stem(p)}} p_n$.
Define in $\V[\co]$,
$$q_1=\bigcup_{n\in A} p_n \quad \text{and}\quad  q_2=\left\{v \in p^\star: \lft1(\co \rest
|v|,v\rgt1) \text{ is good}\right\},$$
and let $q = q_1\cap q_2$.
Since the nodes corresponding to the elements of $A$ were good, it follows
that $q \in \Miller$ and $q \geq_0 p$.
In addition,
$q \forces_{\Miller} \forall n \ \co(n)\leq \dot{\mathbf m}(n)$ (by the
choice of $q_2$) and
$q \forces_{\Miller} \co \not\in \dot{F}$ (by the choice of $q_1$).
Since $\co$ was arbitrary, the proof is finished.
\end{proof}

Finally we show:
\begin{lemma}\label{three}
Suppose that $p \in \Miller$, $s \in \Cohen$,
  and $\lft1(s\rest |\stem(p)|,\stem(p)\rgt1)$ is good.
Let    $\<\dot{F}_n: n \in \omega\>$ be an
  $\Miller$-name for a sequence of closed nowhere dense subsets of $P_{\dot{\mathbf
      m}}$. 
There exists a $\Cohen$-name $\dot{q}$ for an element of $\Miller$
such that 
\begin{enumerate}
\item $s \forces_{\Cohen} \dot{q} \geq p$,
\item $(s, \dot{q}) \forces_{\Cohen \star \dot{\Miller}} \dot{\mathbf c} \in
  P_{\dot{\mathbf m}},$
\item $(s, \dot{q}) \forces_{\Cohen \star \dot{\Miller}} \dot{\mathbf c}
  \not\in \bigcup_n \dot{F}_n$.
\end{enumerate}
\end{lemma}
\begin{proof}
  The proof is a refinement of the proof of the previous lemma.
Suppose that $\<\dot{F}_n: n \in \omega\>$ is an
  $\Miller$-name for a sequence of closed nowhere dense subsets of $P_{\dot{\mathbf
      m}}$. Without loss of generality we can assume that
  $\forces_{\Miller} \forall n \ \dot{F}_n \subseteq \dot{F}_{n+1}$.
Find an 
$\Miller$-name  $\<\dot{s}_n: n \in \omega\>$ such that 
\begin{enumerate}
\item $p \forces_{\Miller} \forall v \in \prod_{j<n} \dot{\mathbf
    m}(j) \ [v^\frown \dot{s}_n] \cap \dot{F}_n=\emptyset$,
\item $p \forces_{\Miller} \forall n \ \lft1(\dot{s}_n, \dot{\mathbf m}
  \rest [n, n+|\dot{s}_n|\rgt1) \text{ is good}$.
\end{enumerate}
Build by induction a sequence of conditions $\<p_n: n \in \omega\>$
such that 
\begin{enumerate}
\item $p_0=p$,
\item $p_{n+1} \geq_{n} p_n$,
\item if $v \in \spli_n(p_{n+1})$ and $k \in U^{p_{n+1}}_v$ then there exists
  a sequence $s^{v^k}_k\in \omega^{<\omega}$ such that 
  \begin{enumerate}
  \item $(p_{n+1})_{v^k} \forces_{\Miller} \dot{s}_k=s^{v^k}_k$,
  \item $|v^k| \geq k+\left|s^{v^k}_k\right|$.
  \end{enumerate}
\end{enumerate}
As in the previous lemma, it follows that for $v \in \spli_n(p_{n+1})$
and $k \in U^{p_{n+1}}_v$,
$\left(s^{v^k}_k, v^k \rest \left[k, k+\lft1|s^{v^k}_k\rgt1|\right)\right)$
is good.

The construction is straightforward; the first step is essentially described in the
previous lemma.
Let $p^\star=\bigcap_n p_n$ and
let $\co \supseteq s$ be a Cohen  real over $\V$.
Working in $\V[\co]$ define for each $n \in \omega $ and  $v \in \spli_n(p^\star)$:
$$A^v=\left\{k \in U^{p^\star}_v \setminus n: \lft1(\co \rest |v^k|, v^k\rgt1) \text{ is
  good and } \co \rest \lft1[k, k+|s^{v^k}_k|\rgt1)=s^{v^k}_k\right\}.$$

As before, it follows that $A^v$ is infinite in $\V[\co]$ for every $v
\in \spli(p^\star)$.
Finally,
let $q \geq p^\star \geq p$ be defined so that
for every $v \in \spli(q)$, $U^q_v=A^v$.

It follows from the definition of $q$ that $\V[\co] \thinks
q \forces_\Miller \co \in P_{\dot{\mathbf m}}$. On the other hand, for
every $v \in \spli_n(q)$, $\V[\co] \thinks q_v \forces_\Miller \exists
k > n \ \co \not\in \dot{F}_k$. Thus 
$$\V[\co] \thinks q \forces_\Miller \exists^\infty n \
\co\not\in \dot{F}_n.$$
Since the sets $\dot{F}_n$ are
 increasing, we conclude that
$\V[\co] \thinks q \forces_\Miller \co \not\in \bigcup_n \dot{F}_n.$
\end{proof}
Now we are ready to prove theorem \ref{crucial}.
Suppose that $X \in \V$, $X \subseteq \omega^\omega $ is not
meager. 
We will show that
$\forces_{\Miller} X \cap P_{\dot{\mathbf m}} \text{ is not meager in }
  P_{\dot{\mathbf m}}$.
Suppose otherwise and let $\<\dot{F}_n: n \in \omega\>$ be an
$\Miller$-name for a meager  set in $P_{\dot{\mathbf m}}$  such that for some
$p \in \Miller$, 
$$p\forces_{\Miller} X \cap P_{\dot{\mathbf m}} \subseteq \bigcup_{n
  \in \omega} \dot{F}_n.$$

Let $N \prec {\mathbf H}(\chi)$ be a countable
elementary submodel  containing $p$, $X$, $\<\dot{F}_n:n \in \omega\>$, etc.
Since $X$ is not meager there exists a real $\co \in X$ which is Cohen
over  $N$. 
By lemma \ref{three} there exists a condition $q \geq p$, $q \in N[\co]$
such that 
$$q \forces_{\Miller} \co \in P_{\dot{\mathbf m}} \setminus \bigcup_{n
  \in \omega} \dot{F}_n.$$
This contradicts the choice of $p$ and finishes the proof.
\end{proof}

Note that in fact we have showed the following:
\begin{theorem}
  Suppose that $\<{\mathcal P}_\alpha, \dot{{\mathcal Q}}_\alpha:
  \alpha<\omega_2\>$ is a countable support iteration of proper
  forcing notions such that:
  \begin{enumerate}
  \item $\{\alpha<\omega_2: \forces_\alpha \dot{{\mathcal Q}}_\alpha \simeq
    \Miller\}$ is cofinal in $ {\boldsymbol\aleph}_2$,
  \item $\forces_\alpha \dot{{\mathcal Q}}_\alpha \text{ preserves
        nonmeager sets}$ for $ \alpha < \omega_2$.
  \end{enumerate}
Then $\V^{{\mathcal P}_{\omega_2}} \thinks \text{Axiom }\Pp$.
\end{theorem}

{\bf Acknowledgements}:
The work was done while I  was spending  sabbatical year at 
 the 
Rutgers University and the College of Staten
Island, CUNY, and I  thank their mathematics departments
for the support. 
I am also grateful to  Andrzej Ros{\l}anowski for careful proofreading  early
versions of this paper and giving many valuable suggestions. 

%\bibliographystyle{plain}
%\bibliography{/home/tomek/biblio/bibliography.bib}

\begin{thebibliography}{10}

\bibitem{BJbook}
Tomek Bartoszy\'{n}ski and Haim Judah.
\newblock {\em Set {T}heory: on the structure of the real line}.
\newblock A.K. Peters, 1995.

\bibitem{JuSh478}
Haim Judah and Saharon Shelah.
\newblock Killing {L}uzin and {S}ierpi\'{n}ski sets.
\newblock {\em Proceedings of the American Mathematical Society}, 120:917--920,
  1994.

\bibitem{Kur66Top}
K.~Kuratowski.
\newblock {\em Topology}, volume~I.
\newblock London; Panstwowe Wydawnictwo Naukowe, Warsaw, 1966.

\bibitem{Mar68}
Edward Marczewski~(Szpilrajn).
\newblock Probleme 68.
\newblock {\em Fundamenta Mathematicae}, 25:579, 1935.

\bibitem{Mil84Spe}
Arnold~W. Miller.
\newblock Special subsets of the real line.
\newblock In K.~Kunen and J.~E. Vaughan, editors, {\em Handbook of Set
  Theoretic Topology}, pages 201--235. North-Holland, Amsterdam, 1984.

\bibitem{Paw89Pro}
Janusz Pawlikowski.
\newblock Products of perfectly meager sets and {L}usin's function.
\newblock {\em Proceedings of the American Mathematical Society},
  107(3):811--815, 1989.

\bibitem{Rec91Pro}
Ireneusz Rec{\l}aw.
\newblock Products of perfectly meagre sets.
\newblock {\em Proceedings of the American Mathematical Society},
  112(4):1029--1031, 1991.

\bibitem{Sie}
Waclaw Sierpi\'nski.
\newblock Sur une probleme de {M}.{K}uratowski concernant la propriete de
  {B}aire des ensambles.
\newblock {\em Fundamenta Mathematicae}, 22:262--266, 1934.

\bibitem{Srivastava}
S.M Srivastava.
\newblock {\em A Course on Borel sets}.
\newblock Graduate Texts in Mathematics. Springer Verlag, 1998.

\bibitem{Zakr}
Piotr Zakrzewski.
\newblock Universally meager sets.
\newblock to appear in Proceedings of the American Mathematical Society.

\end{thebibliography}

\end{document}